\documentclass[12pt]{amsart}

\usepackage{amsmath,amssymb}
\usepackage{amsfonts}
\usepackage{amsthm}
\usepackage{latexsym}
\usepackage{graphicx}
\usepackage{epsfig}

\newtheorem{theorem}{Theorem}[section]
\newtheorem{lemma}[theorem]{Lemma}

\makeatletter \@addtoreset{equation}{section} \makeatother

\theoremstyle{remark}
\newtheorem{remark}{Remark}[section]

\def\ddt{\frac{d}{dt}}

\begin{document}

\title{Backward Ricci Flow on Locally Homogeneous Three-manifolds}

\author{Xiaodong Cao$^*$}
\thanks{$^*$Research
partially supported by the Jeffrey Sean Lehman Fund from Cornell
University}

\address{Department of Mathematics,
  Cornell University, Ithaca, NY 14853}
\email{cao@math.cornell.edu, lsc@math.cornell.edu}


\author{Laurent Saloff-Coste$^{\flat}$}
\thanks{$^{\flat}$Research
partially supported by NSF grant DMS 0603866}


\renewcommand{\subjclassname}{%
  \textup{2000} Mathematics Subject Classification}
\subjclass[2000]{Primary 53C44}

\date{Oct. 15th,  2008}

\maketitle

\markboth{Xiaodong Cao,  Laurent Saloff-Coste} {Backward Ricci
Flow}

\begin{abstract}  In this paper, we study the backward Ricci
 flow on locally homogeneous $3$-manifolds. We describe
the long time behavior and show that, typically and after a proper
re-scaling, there is convergence to a sub-Riemannian geometry. A
similar behavior was observed by the authors in the case of the
cross curvature flow.
\end{abstract}

\section{\bf Introduction}
\subsection{The Ricci flow}
In \cite{ij92}, J. Isenberg and M. Jackson studied the Ricci flow
on homogeneous $3$-manifolds. As homogeneous $3$-manifolds are the
models and building blocks of the geometrization of $3$-manifolds,
it is natural and important to study the behavior of various
geometric flows in this basic case. See \cite[Chapter 1]{ck04}.
Further studies are in \cite{km2001}, \cite{Lott} and
\cite{gli08}.

For obvious reasons, works have focussed on the forward behavior
of the Ricci flow although, in the homogeneous case, the flow
reduces to an ODE system and there is no obstruction to the study
of the backward flow. In \cite{cnsc1} and \cite{csc2}, the authors
studied the forward and backward limits of the cross curvature
flow on  homogeneous $3$-manifolds. Indeed, in the case of the
cross curvature flow it is not entirely clear which direction is
more natural. The results obtained in \cite{csc2} suggest that the
backward behavior of the Ricci flow should be studied as well and
this is the subject of this paper.\\

Recall that the Ricci flow on a manifold is a flow of Riemannian
metric $g(t)$ satisfying the equation
$$\frac{\partial g}{\partial t}= - 2 Rc ,\;\; g(0)=g_0,$$
where $Rc$ denotes the Ricci curvature tensor (in this instance,
the Ricci curvature tensor of the metric $g(t)$). This can be
normalized in various ways by setting
$\widetilde{g}(\widetilde{t})=\psi(t)g(t)$,
$\widetilde{t}=\int_0^t\psi(s)ds$. Setting
$\widetilde{\psi}(\widetilde{t})=\psi(t)$, we obtain
$$\frac{\partial \widetilde{g}}{\partial \widetilde{t}}=
- 2 \widetilde{Rc} +\left(\frac{\partial\ln
\widetilde{\psi}}{\partial \widetilde{t}}\right) \widetilde{g}
,\;\; \widetilde{g}(0)=g_0.$$ For compact manifolds, the customary
normalization uses $\frac{1}{\psi}\frac{\partial\psi}{\partial
t}=\frac{2r}{3}$, where $r$ is the average of the scalar curvature
$R$, in which case,
$\frac{1}{\widetilde{\psi}}\frac{\partial\widetilde{\psi}}{\partial
\widetilde{t}}=\frac{2\widetilde{r}}{3}$. This normalization keeps
the volume constant under the flow. In the case of the locally
homogeneous manifolds, we can use this normalization even in the
non-compact case since the scalar curvature is constant. Hence,
following \cite{ij92}, we will study the flow
\begin{equation}\label{normRF} \frac{\partial g}{\partial t}= - 2
Rc +\frac{2}{3}R g ,\;\; g(0)=g_0.\end{equation}

\subsection{The backward behavior of the Ricci flow}
There are $9$ types of locally homogeneous $3$-manifolds and these
are split into two families. The first family contains the
manifolds covered by the hyperbolic $3$-space $H_3$, and the
product geometries of type $H_2\times \mathbb R$ and $\mathbb
S_2\times \mathbb R$. The second family corresponds to those
geometries whose universal cover is a group itself. They are:
$\mathbb R^3$, $\mbox{SU}(2,\mathbb R)$;
$\widetilde{\mbox{SL}(2,\mathbb R)}$; $E(1,1)= \mbox{Sol}$, i.e.,
the group of isometries of a flat Lorentz plane;
$\widetilde{E(2)}$, the universal cover of group of isometries of
the plane; the Heisenberg group. This second family is referred to
as the Bianchi case (see \cite{ij92}). In the Bianchi case, given
a metric $g_0$, Milnor \cite{milnor76} provides a frame
$(f_1,f_2,f_3)$ in which both the metric  and the Ricci tensors
are diagonalized. As this property is preserved by the Ricci flow,
writing
$$g=A f^1\otimes f^1+Bf^2\otimes f^2 +Cf^3\otimes f^3,$$
the Ricci flow becomes an ODE system in $(A,B,C)$. Furthermore,
Milnor's paper \cite{milnor76}  provides the computation of the
Ricci tensor in each case so that the ODE system in question can
be written down explicitly. The simplest non-trivial case is the
Heisenberg group. Given a metric $g_0$ on the Heisenberg group (or
on a $3$-manifold of Heisenberg  type), we fix a Milnor frame
$\{f_i\}_1^3$ such that $[f_2,f_3]=2f_1,
\;\;[f_3,f_1]=0,\;\;[f_1,f_2]=0.$ Using \cite{milnor76}, the ODE
system for the normalized Ricci flow is given by
\begin{equation}\label{nil+RF}
\left\{
\begin{aligned}
\frac{dA}{dt} =& -\frac{16}{3}\frac{A^3}{A_0B_0C_0}, \\
\frac{dB}{dt} =&+\frac{8}{3}\frac{A^2B}{A_0B_0C_0},\\
\frac{dC}{dt} =& +\frac{8}{3}\frac{A^2C}{A_0B_0C_0},
\end{aligned}
\right .
\end{equation}
where we used the fact that, under (\ref{normRF}),
$ABC=A_0B_0C_0$. Let $R_0= -\frac{2A_0}{B_0C_0}<0$ be the initial
scalar curvature.
 Then (\ref{nil+RF}) admits a completely
explicit maximal solution defined on $(3/(16 R_0),+\infty)$ and
given by
\begin{equation}\nonumber
\left\{
\begin{aligned}
A(t) =& A_0 (1-(16/3) R_0 t)^{-1/2}, \\
B(t) =& B_0 (1- (16/3) R_0 t)^{1/4},\\
C(t) =& C_0 (1-(16/3) R_0 t)^{1/4}.
\end{aligned}
\right .
\end{equation}

Observe that, when $t$ tends to $3/(16R_0)=-T_b$, the metric
$\overline{g}(t)=(C_0)/C(t))g(t)$ converges to
$$\infty f^1\otimes f^1+B_0f^2\otimes f^2 +C_0f^3\otimes f^3,$$
which can be interpreted as describing a sub-Riemannian geometry
on the Heisenberg group.  The point of this paper is to show that
this behavior in the backward direction is typical  for all
locally homogeneous manifolds corresponding to the Bianchi cases
described above except for those corresponding to the trivial case
$\mathbb R^3$.

\begin{theorem}\label{th-main}
Let $(M,g_0)$ be a locally homogeneous $3$-manifold with universal
cover $\mbox{\em SU}(2,\mathbb R)$, $\widetilde{\mbox{\em
SL}(2,\mathbb R)}$, $E(1,1)= \mbox{\em Sol}$, $\widetilde{E(2)}$,
or the Heisenberg group. Let $g(t)$, $(-T_b,T_f)$ be a maximal
solution of the normalized Ricci flow {\em (\ref{normRF})}. Let
$d(t)$ be the associated distance function on $M$. Assume $g_0$ is
generic among all locally homogeneous metric on $M$. Then $T_b$ is
finite and there exists a function $r(t):(-T_b,0]\rightarrow
(0,\infty)$ such that, as $t$ tends to $-T_b$, the metric spaces
$(M,r(t)d(t))$ converge uniformly to a sub-Riemannian space
$(M,d_b)$ whose tangent cone at any point is the Heisenberg group
equipped with its natural sub-Riemannian metric.
\end{theorem}

By definition, the uniform convergence of metric spaces $(M, d_t)$
to $(M,d)$ means the uniform convergence over compact sets of
$(x,y) \rightarrow d_t (x,y)$ to $(x,y) \rightarrow d (x,y)$.

The present paper proves this theorem in all cases except
$\widetilde{\mbox{ SL}(2,\mathbb R)}$. For manifolds covered by
$\widetilde{\mbox{ SL}(2,\mathbb R)}$, we prove the result under
the additional assumption that there exists a time $t_0$ such that
either $A(t_0) \ge \max\{B(t_0), C(t_0)\}$ or $A(t_0) \le
|B(t_0)-C(t_0)|$. In the paper \cite{cgsc08}, we show that this
condition is always satisfied except for a hypersurface of initial
conditions.

The proof of this theorem proceeds by inspection of the different
cases. It would be more elegant to have an argument covering all
cases at once. However, the existence of exceptional sets of
initial conditions for which the general result fails indicates
that it is unlikely that such treatment is possible. Indeed, the
exceptional sets of initial conditions are very much
case-dependent, see the more precise statements in the different
sections below.

The results obtained in each of the different cases are more
precise than stated in Theorem \ref{th-main}. They describe the
asymptotic behavior of each of the metric components in a fixed
Milnor frame. This is useful in exploring the Ricci flow on
homogeneous $3$-manifolds under more sophisticated scaling
procedures. See \cite{km2001,Lott,gli08}.

Together, the study of the forward normalized Ricci flow (see
\cite{ij92, Lott, km2001}) and this paper, give a description of
the asymptotic behaviors of the Ricci flow on homogeneous
$3$-manifolds for both the forward and backward directions. For
instance, the solutions of the forward normalized Ricci flow
always exist for all (positive) time in the Bianchi classes
(\cite{ij92}).

\subsection{Sub-Riemannian geometries}
Our main result, Theorem \ref{th-main}, refers to the notion of
sub-Riemannian geometry, a term that we now explain in the present
context. The typical behavior (possibly after some re-scaling) of
the evolving metric
$$g=A f^1\otimes f^1+Bf^2\otimes f^2
+Cf^3\otimes f^3,$$ at the end points of a maximal existence
interval is that some of the coefficients $A,B,C$ either vanish or
tend to $\infty$. When a coefficient vanishes and the manifold is
compact, the phenomenon can be interpreted as a dimensional
collapse. Naively, at least one direction disappears. To interpret
the situation when a coefficient tends to infinity, it is useful
to look at the dual tensor
$$Q=A^{-1} f_1\otimes f_1+B^{-1}f_2\otimes f_2
+C^{-1}f_3\otimes f_3$$ defined on the co-tangent bundle. Suppose
that $A$ tends to infinity whereas $B,C$ have finite limits
$B_*,C_*$. Then the tensor $Q$ tends to
$$Q_*=
B_*^{-1}f_2\otimes f_2 +C_*^{-1}f_3\otimes f_3.$$ If it turns out
that $[f_2,f_3]=2\epsilon_1 f_1$ with $\epsilon_1\neq 0$, then the
tensor $Q_*$ induces a natural distance function $d_*$ on the
underlying manifold $M$. This distance can be computed by
minimizing the length of the so-called horizontal curves, i.e.,
those curves that stay tangent to the linear span of $f_2,f_3$.
The associated ``geometry'' is called a sub-Riemannian geometry.
See \cite{mon02} for a book length introduction to sub-Riemannian
geometry and \cite{csc2} for some details directly relevant to the
present situation. Let us note here that the convergence
$Q\rightarrow Q_*$ translates quite easily in the uniform
convergence over compact sets on $M\times M$ of the associated
distance functions. This explains the conclusion of Theorem
\ref{th-main}.

\subsection{The normalized backward Ricci flow}
In order to study the backward behavior of the Ricci flow, it is
convenient to reverse time and consider the solution of the
positive normalized Ricci flow equation
\begin{equation}\label{norm+RF} \frac{\partial g}{\partial t}=  2
Rc -\frac{2}{3}R g ,\;\; g(0)=g_0.\end{equation} We let $T_+\in
[0,+\infty]$ be the maximal existence time for this equation. The
rest of this paper is devoted to the asymptotic properties of this
flow when $t\rightarrow T_+$ in the case of $\mbox{SU}(2)$,
$\mbox{E}(1,1)$, $\widetilde{\mbox{E}(2)} $ and
$\widetilde{\mbox{SL}(2,\mathbb R)}$. This includes determining
whether $T_+$ is finite or infinite. The results are stated
explicitly for the flow on each of these groups but, in each case,
it holds in the same form on any locally homogeneous $3$-manifold
covered by  the corresponding group. In each case, we write the
solution of (\ref{norm+RF}) in the form $$g=A f^1\otimes
f^1+Bf^2\otimes f^2 +Cf^3\otimes f^3$$ in a Milnor frame
$(f_1,f_2,f_3)$ for $g_0$. Under (\ref{norm+RF}), $ABC=A_0B_0C_0$
is a constant. In the rest of this paper, we assume the
normalization $A_0B_0C_0=4$. This choice is made so that the ODE
systems are the same as in \cite{ij92}, despite the fact that the
frame we use here have a different normalization than those used
in \cite{ij92}.

If $A, B, C$ is the solution under $A_0B_0C_0=4$ and
$\widetilde{A}, \widetilde{B}, \widetilde{C}$ is the solution with
$\widetilde{A_0}=\lambda A_0, \widetilde{B_0}=\lambda B_0,
\widetilde{C_0}=\lambda C_0$, then $\widetilde{A}(t)=\lambda A(t/
\lambda)$, $\widetilde{B}(t)=\lambda B(t/ \lambda)$ and
$\widetilde{C}(t)=\lambda C(t/ \lambda)$.

\section{\bf The  normalized positive Ricci flow on $\mbox{SU}(2)$}

Given a metric $g_0$ on $\mbox{SU}(2)$, we fix a Milnor frame such
that $[f_i,f_j]=2f_k$ for all cyclic permutations of the indices.
This section is devoted to the proof of the following result.
\begin{theorem} \label{thm-su2} Let $g_0$ be an homogeneous metric
on $\mbox{\em SU}(2)$ with associated Milnor frame $(f_1,f_2,f_3)$
and $g_0=A_0f^1\otimes f^1+B_0 f^2\otimes f^2+C_0 f^3\otimes f^3$
with $A_0B_0C_0=4$. Let $g(t)=A(t)f^1\otimes f^1+B(t) f^2\otimes
f^2+C(t) f^3\otimes f^3 $, $t\in [0,T_+)$ be the maximal forward
solution of the positive normalized Ricci flow {\em
(\ref{norm+RF})} with $g(0)=g_0$. Assume that $A_0\ge B_0\ge C_0$.
\begin{enumerate}
\item If $A_0=B_0=C_0$ then $T_+=\infty$ and $g(t)=g_0$, $t\in
[0,\infty)$. \item If $A_0=B_0>C_0$ then $T_+=\infty$, $A=B>C$
and, as $t$ tends to infinity, $A\sim \frac{8}{3}t$, $C\sim
\frac{9}{16}  t^{-2}.$ \item If $A_0>B_0\ge C_0$ then $T_+$ is
finite, $A>B\ge C$ and there are constants $\eta_1,\eta_2\in
(0,\infty)$ such that
$$A\sim \frac{\sqrt{6}}{4}(T_+-t)^{-1/2},\;\; B\sim
\eta_1(T_+-t)^{1/4},\;\; C\sim \eta_2(T_+-t)^{1/4}$$ as $t$ tends
to $T_+$.
\end{enumerate}
Let $d(t)$ be the distance function associated to
$(B_0/B(t))g(t)$. In case {\em (3)}, the metric space $(\mbox{\em
SU}(2),d(t))$ converges uniformly as $t\rightarrow T_+$  towards
the sub-Riemanninan metric space $(\mbox{\em SU}(2),d_*)$ where
$d_*$ is the sub-Riemannian distance associated with
$$Q_*=B^{-1}_0 f_2\otimes f_2+
\eta_1\eta^{-1}_2B_0^{-1}f_3\otimes f_3.$$
\end{theorem}

\begin{remark}
Consider a maximal solution $$g_f(t)=A(t)f^1\otimes f^1+B(t)
f^2\otimes f^2+C(t) f^3\otimes f^3, \;\; t\in (-T_+,\infty)$$ of
the forward normalized Ricci flow {\mbox (\ref{normRF})}. Let
$\overline{g} (t)=(B_0/B(t))g_f(t)$. Isenberg and Jackson
\cite{ij92} shows that $A-C\leq (A_0-C_0)e^{-2C_0^2t}, \forall
t\geq 0$, if $A_0\ge B_0\ge C_0$ (this order is preserved by the
flow). Hence, in the forward direction, $\overline{g} (t)$
converges exponentially fast to the round metric whereas Theorem
\ref{thm-su2} describes the backward behavior. In the generic case
$A_0>B_0\ge C_0$, $\overline{g} (t)$ converges to a sub-Riemannian
metric as $t\rightarrow -T_+$.
\end{remark}

The sectional curvatures are (see, e.g., \cite [pg. 12]{ck04}
\begin{align*}
K(f_2 \wedge f_3)&=\frac{(B-C)^2}{ABC}-\frac{3A}{BC}+\frac2B+\frac2C,\\
K(f_3 \wedge f_1)&=\frac{(C-A)^2}{ABC}-\frac{3B}{CA}+\frac2A+\frac2C,\\
K(f_1 \wedge
f_2)&=\frac{(A-B)^2}{ABC}-\frac{3C}{AB}+\frac2A+\frac2B.
\end{align*}

From the sectional curvatures given above, we easily obtain the
ODEs corresponding to the flow, under the normalization $ABC=4$,
namely,
\begin{equation}\label{pdesu2}
\left \{
\begin{aligned}
\frac{dA}{dt}=&-\frac{2}{3}A[-A(2A-B-C)+(B-C)^2],\\
\frac{dB}{dt}=&-\frac{2}{3}B[-B(2B-A-C)+(A-C)^2],\\
\frac{dC}{dt}=&-\frac{2}{3}C[-C(2C-A-B)+(A-B)^2].
\end{aligned}
\right .
\end{equation}

Without loss of generality we may assume that $A_0\ge B_0\ge C_0$.
As
\begin{align}
\ddt (A-C)=&\frac{2}{3}(A-C)[2A^2+2AC+2C^2-(A+B+C)B],
\label{a-c}\\
\ddt (A-B)=&\frac{2}{3}(A-B)[2A^2+2AB+2B^2-(A+B+C)C],
\label{a-b}\\
\ddt (B-C)=&\frac{2}{3}(B-C)[2B^2+2BC+2C^2-(A+B+C)A], \label{b-c}
\end{align}
it is easy to see that $A\ge B\ge C$ is preserved along the flow.
This yields the following lemma.
\begin{lemma}\label{lem-mon}
Assume that $A_0\ge B_0\ge C_0$. Then $A$, $A-B$ and $A-C$ are all
nondecreasing along the flow and $C$ is non-increasing.
\end{lemma}

We now consider three cases. The first case is when $A_0=B_0=C_0$.
Then $A(t)=B(t)=C(t)=A_0$ and the solution exists for all time.\\

The second case is when  $A_0=B_0>C_0$. Then $A(t)=B(t)$ as long
as the solution exists and we have
\begin{equation}\label{Sua=b}
\left \{
\begin{aligned}
\frac{dA}{dt}=&\frac{2}{3}AC(A-C),\\
\frac{dC}{dt}=&-\frac{4}{3}C^2(A-C).
\end{aligned}
\right .
\end{equation}
In this case, $A$ is increasing, $C$ is decreasing and
$A^2C=A_0^2C_0$ along the flow.
\begin{lemma}
If $A_0=B_0>C_0$, then $T_+=\infty$,
 $A\sim \frac83  t$, and $
C\sim \frac{9}{16}  t^{-2}$ as $t$ tends to infinity.
\end{lemma}
\begin{proof}
Since $\frac{dA}{dt}=\frac{2}{3}A_0^2C_0-\frac{2}{3}AC^2,$  if
$T_+<\infty$, then $\lim_{T_+} A<\infty$, and $\lim_{T_+} C>0$.
This contradicts the assumption that $T_+$ is the maximal
existence time. Hence $T_+=\infty$. As $A$ is increasing, $C$
decreasing and $A^2C$ constant, it follows from (\ref{Sua=b}) that
$\lim_{\infty} A= \infty$, and thus $\lim_{\infty} C=0$. Moreover,
$\lim_{\infty} AC^2= 0$. Now the asymptotic for $A$ and $C$
follows from (\ref{Sua=b}) which yields $\ddt A \sim \frac23
A_0^2C_0$, and $A^2C=A_0^2C_0=4$.
\end{proof}

We now focus on the third case, the generic case.
\begin{lemma}\label{case31}
Assume that $A_0> B_0\ge C_0$. Then $T_{+}<\infty$.
\end{lemma}
\begin{proof}
Assume that $T_{+}=\infty$. We have $$\ddt A=\frac23A[A(A-B)
+A(A-C)-(B-C)^2]\geq \frac23A^2(A-B),$$ and $$\ddt C
<-\frac23C(A-B)^2.$$ Since, by Lemma \ref{lem-mon}, both $A$ and
$A-B$ are nondecreasing, it follows that $\lim_{\infty} A=\infty$,
$\lim_{\infty} C=0$. Now, (\ref{a-b}) implies  that
$$\ddt \ln (A-B)\geq \frac{2}{3}(2A^2+B^2),$$ hence
$\lim_{\infty} (A-B)=\infty$. Since
\begin{equation}\label{lnb}
\ddt \ln B
=-\frac23[B(A-B)+(A-B)^2+2(A-B)(B-C)-C(B-C)],
\end{equation}
 this shows that $B$ is non-increasing for $t$ large enough,
 hence bounded. So, we have
\begin{equation*}
\frac{dA}{dt}\sim \frac{4}{3}A^3.
\end{equation*}
But this shows that there exists a finite time
$T_0$, such that $\lim_{T_0} A=\infty$, this
contradicts our assumption that $T_{+}=\infty$.
\end{proof}

\begin{lemma}\label{}
Assume  $A_0> B_0\ge C_0$. Then $\lim_{T_{+}} A=\infty$,
$\lim_{T_{+}} B=\lim_{T_{+}} C=0$.
\end{lemma}
\begin{proof}
Assume that $\lim_{T_{+}} C >0$. As $A> B\geq C$ and that $T_{+}$
is finite, we must have $\lim_{T_{+}} A=\infty$.  We have
$$\ddt \ln A <\frac43 A^2 \le \ddt \ln (A-B).$$
It follows that $\lim_{T_{+}} (A-B)=\infty$. By (\ref{lnb}), $B$
is non-increasing for $t$ close to $T_{+}$ and hence bounded from
above. This shows that $\ddt A \sim \frac43A^3$ and  thus that
$A^{-2} \sim \frac83 (T_{+}-t)$.
 Hence $$\ddt \ln C \sim
-\frac14(T_{+}-t)^{-1}.$$ This contradicts  $\lim_{T_{+}} C
>0$ and we conclude that $\lim_{T_{+}} C=0$.

Now by (\ref{lnb}) we can see that $B$ is bounded from above. So,
if $\lim_{T_{+}} A<\infty$ then $\ddt C\sim -\eta C,$ for some
constant $\eta\in (0,\infty)$. This contradicts $\lim_{T_{+}}
C=0$. So we conclude that $\lim_{T_{+}} A=\infty$.

To show that $\lim_{T_{+}} B=0$, notice that (\ref{lnb}) implies
that $B$ is non-increasing for $t$ close enough to $T_+$. As
$$\ddt \ln (AB^2)=2(B-C)(B+C-A),$$ we obtain that
$AB^2$ is bounded from above
on $[0,T_+)$, hence $\lim_{T_{+}} B=0$.
\end{proof}

\begin{lemma}Assume that $A_0> B_0\ge C_0$, then there exist
$\eta_1,~ \eta_2 \in (0,\infty)$, such that
$$ A \sim \frac{\sqrt{6}}{4}(T_{+}-t)^{-1/2},\;\;
B\sim \eta_1 (T_+-t)^{1/4},\;\; C\sim \eta_2 (T_+-t)^{1/4}.$$
\end{lemma}
\begin{proof}
The first statement follows directly from
$$\ddt A \sim \frac43A^3.$$
To obtain the asymptotic behavior for $B, C$, notice that
$$\ddt \ln (AC^2)=2(B-C)(A-B-C),\;\;
\ddt \ln (AB^2)=2(B-C)(B+C-A).
$$
Since  $\lim_{T_{+}} A=\infty$ and $\lim_{T_{+}} B= \lim_{T_{+}}
C=0$, we have $\lim_{T_{+}} (A-B-C)=\infty$. Hence, the equations
above imply that $AB^2$ is non-increasing   and $AC^2$ is
non-decreasing for $t$ close to $T_+$.  But $B\ge C$, so
$$0<\lim_{T_{+}} AC^2\le \lim_{T_{+}} AB^2<\infty.$$
It follows that $\lim_{T_{+}} B\sim \eta_1 (T_+-t)^{1/4}$ and
$\lim_{T_{+}} C\sim \eta_2 (T_+-t)^{1/4}$.
\end{proof}

This finishes the proof of Theorem \ref{thm-su2}.

\section{\bf The  normalized positive Ricci flow on $\mbox{E}(1,1)$
(Sol geometry)}

Given a metric $g_0$ on $\mbox{E}(1,1)$, we fix a Milnor frame
such that $[f_2,f_3]=2f_1$, $[f_3,f_1]=0$, $[f_1,f_2]=-2f_3$. This
section is devoted to the proof of the following result.
\begin{theorem}\label{thE11}
Let $g_0$ be an homogeneous metric on $\mbox{\em E}(1,1)$ with
associated Milnor frame $(f_1,f_2,f_3)$ and $g_0=A_0f^1\otimes
f^1+B_0 f^2\otimes f^2+C_0 f^3\otimes f^3$ with $A_0B_0C_0=4$. Let
$g(t)=A(t)f^1\otimes f^1+B(t) f^2\otimes f^2+C(t) f^3\otimes f^3
$, $t\in [0,T_+)$ be the maximal forward solution of the positive
normalized Ricci flow {\em (\ref{norm+RF})} with $g(0)=g_0$.
Assume that $A_0\ge C_0$.
\begin{enumerate}
\item If $A_0=C_0$ then $T_+=\frac{3}{32} B_0$ and
$$A(t)=C(t)=
\frac{\sqrt{6}}{4}(T_+-t)^{-1/2}, \;\;B(t)=\frac{32}{3}(T_+-t),
\;\;t\in [0,T_+).$$ \item If $A_0>C_0$ then $T_+<\infty$ and, as
$t$ tends to $T_+$, there are constants $\eta_1,\eta_2\in
(0,\infty)$ such that
$$A\sim \frac{\sqrt{6}}{4}(T_+-t)^{-1/2},\;\; B\sim
\eta_1(T_+-t)^{1/4},\;\; C\sim \eta_2(T_+-t)^{1/4}.$$
\end{enumerate}
Let $d(t)$ be the distance function associated to
$(B_0/B(t))g(t)$. In case {\em (2)}, the metric space $(\mbox{\em
E}(1,1),d(t))$ converges uniformly as $t\rightarrow T_+$  towards
the sub-Riemanninan metric space $(\mbox{\em E}(1,1),d_*)$ where
$d_*$ is the sub-Riemannian distance associated with
$$Q_*=B^{-1}_0 f_2\otimes f_2+
\eta_1\eta^{-1}_2B_0^{-1}f_3\otimes f_3.$$
\end{theorem}

\begin{remark}
For the forward normalized Ricci flow (\ref{normRF}), Isenberg and
Jackson \cite{ij92} show that the solution exists for all time and
presents a cigar degeneracy.
\end{remark}

The sectional curvatures of $g(t)$ in the frame $(f_i)_1^3$  are:
\begin{eqnarray*}
K(f_2 \wedge f_3) &= & \frac{(A-C)^2-4A^2}{ABC},\\
K(f_3 \wedge f_1) &= & \frac{(A+C)^2}{ABC},\\
K(f_1 \wedge f_2) &= & \frac{(A-C)^2-4C^2}{ABC}.
\end{eqnarray*}
These yield  the equations for the normalized positive Ricci flow
on $E(1,1)$, under the normalization $ABC=4$, namely,
\begin{equation}
\left \{
\begin{aligned}
\frac{dA}{dt} =&\frac{2}{3}A(2A^2+AC-C^2),
\\
\frac{dB}{dt} =&-\frac{2}{3}B(A+C)^2,
\\
\frac{dC}{dt} =&\frac{2}{3}C(2C^2+AC-A^2).
\end{aligned}
\right .
\end{equation}

\begin{lemma}
If $A_0=C_0$, then $T_+=\frac{3}{32}B_0<\infty$. Moreover
$A(t)=C(t)=\frac{\sqrt{6}}{4} (T_+ -t)^{-1/2}$ and
$B(t)=\frac{32}{3} (T_+-t)$, for $t\in [0,T_+)$.
\end{lemma}
\begin{proof}
It is easy to see that $A=C$ as long as the solution exists. As
$\ddt A^{-2}=-\frac83$, we have  $T_+<\infty$ and
$$A=\frac{\sqrt{6}}{4} (T_+ -t)^{-1/2},\;T_+=\frac{3}{8A_0^2}=
\frac{3}{32}B_0.$$ Further, $\ddt B= -\frac{8}{3}BA^2,$ so
$B=\frac{32}{3} (T_+-t)$.
\end{proof}

Without loss of generality, we assume that $A_0>C_0$. This implies
that $A$ is increasing. Note that $B$ is always decreasing.

\begin{lemma}\label{lemE111}
If $A_0>C_0$, then $T_+<\infty$, $\lim_{T_{+}} A=\infty$, and
there exists a time $t_0$ such that $A(t_0)\geq 2C(t_0)$.
\end{lemma}
\begin{proof}
The fact that $T_+<\infty$ follows from $\ddt A>\frac23A^3$. Now
assume that $\lim_{T_{+}} A= A(T_+)<\infty$. Then, since $$\ddt
\ln C= -\frac23 (A+C)(A-2C)>-\frac43 A^2> -\frac43 A(T_+)^2,$$ and
$$\ddt \ln B=-\frac23(A+C)^2>\frac83 A^2> -\frac83 A(T_+)^2,$$
we get that $B\ge \lim_{T_{+}} B=B(T_+)>0$. Similarly, there
exists some constant $\eta>0$ such that $C\in [\eta,A(T_+)]$. This
contradicts the fact that the maximal existence time $T_+$ is
finite. Hence $\lim_{T_{+}} A=\infty$.

To prove the second statement, we assume that $C<A<2C$ for all $t
\in [0,T_+)$. So we have $\lim_{T_{+}} C=\infty$. Since
$$\ddt \ln (A/C)=2(A+C)(A-C)>0,$$ we see that $A/C$ is increasing, so
$A/C >A_0/C_0$ and $A-C=A(1-C/A)>(1-C_0/A_0)A$. Moreover, we have
$$(1-C_0/A_0)\int_0^{T_+} (A+C)A < \int_0^{T_+} (A+C)(A-C) <\frac12 \ln 2.$$ Hence
$$\int_0^{T_+} (A+C)(2C-A)<\int_0^{T_+} (A+C)C<
\int_0^{T_+} (A+C)A <\infty.$$ This contradicts the fact that
$$\ddt \ln C=\frac23(A+C)(2C-A) \mbox{ and }\lim_{T_{+}}
C=\infty.$$ So there exists a time $t_0$ such that $A(t_0)\geq
2C(t_0)$.
\end{proof}

\begin{lemma}Assume that $A_0> C_0$. There exist
$\eta_1,~ \eta_2 \in (0,\infty)$ such that, as $t$ tends to $T_+$,
we have
$$A \sim \frac{\sqrt{6}}{4}(T_{+}-t)^{-1/2},\;\;
B\sim \eta_1 (T_+-t)^{1/4},\;\; C\sim \eta_2 (T_+-t)^{1/4}.$$
\end{lemma}

\begin{proof} By Lemma \ref{lemE111}, there is $t_0$ such that
$A(t_0)\ge 2C(t_0)$. As $\ddt \ln (A/C)=2(A+C)(A-C)>0$, we
conclude that $A(t)\geq 2C(t)$ for $t\in [t_0,T_+)$. Hence $C$ is
non-increasing on $[t_0,T_+)$. As
$$\ddt \ln (AC^3)=\frac23(-A^2+4AC+5C^2),\;\;
\ddt \ln (AB^2)=-2C(A+C),$$ and $\lim_{T_{+}} A=\infty$, it
follows that both $AC^3$ and $AB^2$ are bounded from above, hence
$\lim_{T_{+}} B=\lim_{T_{+}} C=0$.

Next, we show that $\lim_{T_{+}} AB^2=\eta_1$ and $\lim_{T_{+}}
AC^2=\eta_2$. Note that
$$\ddt \ln (AC^2)=2C(A+C) \mbox{ and }  \ddt \ln (AB^2)=-2C(A+C).$$
Hence, it is enough to prove that $\int_0^{T_+} AC<\infty$. As
$\ddt C \sim -\frac23 A^2C,$ and $C>0$, we have
$$\int_0^{T_+} AC<A_0^{-1}\int_0^{T_+} A^2C<\infty.$$
Now, the lemma follows from $\ddt A \sim \frac43A^3$.
\end{proof}
This ends the proof of Theorem \ref{thE11}.

\section{\bf The normalized positive Ricci flow on
$\widetilde{\mbox{E}(2)}$}

Given a left-invariant metric $g_0$ on $\widetilde{\mbox{E}(2)}$,
we fix a Milnor frame $\{f_i\}_1^3$ such that
$$[f_2,f_3]=2f_1, \;\;[f_3,f_1]=2f_2,\;\;[f_1,f_2]=0.$$
The result in this case reads as follows.
\begin{theorem}\label{thE2}
Let $g_0$ be an homogeneous metric on $\widetilde{\mbox{\em
E}(2)}$ with associated Milnor frame $(f_1,f_2,f_3)$ and
$g_0=A_0f^1\otimes f^1+B_0 f^2\otimes f^2+C_0 f^3\otimes f^3$ with
$A_0B_0C_0=4$. Let $g(t)=A(t)f^1\otimes f^1+B(t) f^2\otimes
f^2+C(t) f^3\otimes f^3 $, $t\in [0,T_+)$ be the maximal forward
solution of the positive normalized Ricci flow {\em
(\ref{norm+RF})} with $g(0)=g_0$. Assume that $A_0\ge B_0$.
\begin{enumerate}
\item If $A_0=B_0$ then $T_+=\infty$ and $g(t)=g_0$ on
$[0,\infty)$. \item If $A_0>B_0$ then $T_+<\infty$ and, as $t$
tends to $T_+$, there are constants $\eta_1,\eta_2\in (0,\infty)$
such that
$$A\sim \frac{\sqrt{6}}{4}(T_+-t)^{-1/2},\;\; B\sim
\eta_1(T_+-t)^{1/4},\;\; C\sim \eta_2(T_+-t)^{1/4}.$$
\end{enumerate}
Let $d(t)$ be the distance function associated to
$(B_0/B(t))g(t)$. In case {\em (2)}, the metric space
$(\widetilde{\mbox{\em E}(2)},d(t))$ converges uniformly as
$t\rightarrow T_+$  towards the sub-Riemanninan metric space
$(\widetilde{\mbox{\em E}(2)},d_*)$ where $d_*$ is the
sub-Riemannian distance associated with
$$Q_*=B^{-1}_0 f_2\otimes f_2+
\eta_1\eta^{-1}_2B_0^{-1}f_3\otimes f_3.$$
\end{theorem}

\begin{remark}
Consider a maximal solution $$g_f(t)=A(t)f^1\otimes f^1+B(t)
f^2\otimes f^2+C(t) f^3\otimes f^3, \;\; t\in (-T_+,\infty)$$ of
the forward normalized Ricci flow {\mbox (\ref{normRF})}. Let
$\overline{g} (t)=(B_0/B(t))g_f(t)$. Isenberg and Jackson
\cite{ij92} shows that $A-B\leq (A_0-B_0)e^{-4B_0^2t}, \forall
t\geq 0$, if $A_0\ge B_0$ (this order is preserved by the flow).
Hence, in the forward direction, $\overline{g} (t)$ converges
exponentially fast to the flat metric whereas Theorem \ref{thE2}
describes the backward behavior. In the generic case $A_0>B_0$,
$\overline{g} (t)$ converges to a sub-Riemannian metric as
$t\rightarrow -T_+$.
\end{remark}

In this case, the sectional curvatures are:
\begin{eqnarray*}
K(f_2 \wedge f_3) &=& \frac{1}{ABC} (B-A)(B+3A),\\
 K(f_3 \wedge f_1) &= &
\frac{1}{ABC} (A-B)(A+3B),\\
K(f_1 \wedge f_2) &= & \frac{1}{ABC} (A-B)^2.
\end{eqnarray*}
 Hence the solution
$g(t)=A(t)f^1\otimes f^1+B(t)f^2\otimes f^2+C(t)f^3\otimes f^3$ of
the normalized positive Ricci flow satisfies
\begin{equation}\label{pdee2}
\left\{
\begin{aligned}
\frac{dA}{dt}=&\frac{2}{3}A(2A+B)(A-B),\\
\frac{dB}{dt}=&-\frac{2}{3}B(2B+A)(A-B),\\
\frac{dC}{dt}=&-\frac{2}{3}C(A-B)^2,
\end{aligned}
\right .
\end{equation}
under the normalization $ABC=4$.

 If $A_0=B_0$ we clearly have
$g(t)=g_0$ for all $t\ge 0$. Without loss of generality, we assume
that $A_0>B_0$. Then $A>B$ as long as the solution exists. Hence
$A$ is increasing whereas $B$ and $C$ are decreasing.

\begin{lemma}
If $A_0>B_0$, then $T_+<\infty$, $\lim_{T_+} A=\infty$, $\lim_{T_+} B
=\lim_{T_+} C=0$.
\end{lemma}
\begin{proof}
Since $A-B>0$ is increasing and
$$\ddt (A-B)=\frac43 (A-B)(A^2+AB+B^2)>\frac43 (A-B)^3,$$
we have  $\ddt (A-B)^{-2}<-\frac83$, so $T_+<\infty$.

If $\lim_{T_+} A=A(T_+)<\infty$, then $$\ddt \ln
B>-2A^2>-2A(T_+)^2,$$ and $$\ddt \ln C>-\frac23 A^2>-\frac23
A(T_+)^2.$$ This leads to $\lim_{T_+} B=B(T_+)>0$ and $\lim_{T_+}
C=C(T_+)>0$ and  contradicts the fact that the maximal existence
time $T_+$ is finite.

To prove that $B$ tends to $0$, note that
$$\ddt \ln (AB^2)=-2(A-B)B<0.$$ Hence $AB^2$ is decreasing and
$\lim_{T_+} B=0$. Similarly,
$$\ddt \ln (AC^3)=-\frac23 (A-B)(A-2B)$$ implies  that $AC^3$
is bounded from above. Hence $\lim_{T_+} C=0$.
\end{proof}

\begin{lemma}Assume that $A_0> B_0$. Then there exist
$\eta_1,~ \eta_2 \in (0,\infty)$ such that, as $t$ tends to $T_+$,
$$A \sim \frac{\sqrt{6}}{4}(T_{+}-t)^{-1/2},\;\;
B\sim \eta_1 (T_+-t)^{1/4},\;\; C\sim \eta_2 (T_+-t)^{1/4}.$$
\end{lemma}
\begin{proof}Since $B>0$ and $\ddt B \sim -\frac23 A^2B$, we get
\begin{equation}\label{lim_AB}
\int_0^{T_+} AB< A_0^{-1}\int_0^{T_+} A^2B<\infty.
\end{equation} Observe that $$\ddt \ln
(AB^2)=-2(A-B)B;\;\; \ddt \ln (AC^2)=2 (A-B)B. $$ Hence
(\ref{lim_AB}) implies $$\lim_{T_+} AB^2>0;\;\; \lim_{T_+} AC^2
<\infty.$$  The asymptotic behaviors of $A$, $B$ and $C$ now
follow from $\ddt A \sim \frac43 A^3$.
\end{proof}

This finishes the proof of Theorem \ref{thE2}.

\section{\bf The normalized positive Ricci flow on
$\mbox{SL}(2,\mathbb R)$}

Given a left-invariant  metric $g_0$ on $\mbox{SL}(2,\mathbb R)$,
we fix a Milnor frame $\{f_i\}_1^3$ such that
$$[f_2,f_3]=-2f_1, \;\;[f_3,f_1]=2f_2,\;\;[f_1,f_2]=2f_3$$
and $$g_0=A_0f^1\otimes f^1+B_0f^2\otimes
f^2+C_0f^3\otimes f^3.$$

\begin{theorem}\label{th-sl2r}
Let $g_0$ be an homogeneous metric on $\mbox{SL}(2,\mathbb R)$
with associated Milnor frame $(f_1,f_2,f_3)$ and
$g_0=A_0f^1\otimes f^1+B_0 f^2\otimes f^2+C_0 f^3\otimes f^3$ with
$A_0B_0C_0=4$. Let $g(t)=A(t)f^1\otimes f^1+B(t) f^2\otimes
f^2+C(t) f^3\otimes f^3 $, $t\in [0,T_+)$ be the maximal forward
solution of the positive normalized Ricci flow {\em
(\ref{norm+RF})} with $g(0)=g_0$. Then $T_+<\infty$. Assume that
$B_0\ge C_0$, and set
$$Q=\{(a,b,c)\in \mathbb R^3: a>0,
b\ge c>0\}$$ and
$$\bar{g}(t)=\frac{C_0}{C(t)} g(t).$$
There is a partition of $Q$ into subsets $S_0,Q_1,Q_2$ with
$Q_1,Q_2$ connected such that, as $t$ tends to $T_+$:
\begin{enumerate}
\item  If $(A_0,B_0,C_0)\in Q_1$ then there exist $\eta_1, \eta_2
\in (0,\infty)$ such that $$A\sim
\frac{\sqrt{6}}{4}(T_+-t)^{-1/2},\;\; B\sim
\eta_1(T_+-t)^{1/4},\;\; C\sim \eta_2(T_+-t)^{1/4}.$$ Moreover,
$(M,\overline{g}(t))$ converges uniformly to the sub-Riemannian
metric space $(M, bf_2\otimes f_2+ cf_3\otimes f_3)$ for some
$b,c\in (0,\infty)$.

\item  If $(A_0,B_0,C_0)\in Q_2$ then there exist $\eta_1, \eta_2
\in (0,\infty)$ such that
$$A\sim \eta_1(T_+-t)^{1/4},
\;\; B\sim \frac{\sqrt{6}}{4}(T_+-t)^{-1/2},\;\; C\sim
\eta_2(T_+-t)^{1/4}.$$ Moreover, $(M,\overline{g}(t))$ converges
uniformly to the sub-Riemannian metric space $(M, af_1\otimes f_1+
cf_3\otimes f_3)$ for some $a,c\in (0,\infty)$.

\item If $(A(t),B(t),C(t))\in S_0$ for all $t\in (T_+,0]$ then
$$A\sim \frac{\sqrt{6}}{4}(T_+-t)^{-1/2},
\;\; B\sim \frac{\sqrt{6}}{4}(T_+-t)^{-1/2},\;\; C\sim
\frac{32}{3}(T_+-t).$$
\end{enumerate}
\end{theorem}

\begin{remark} The cases (1)-(2) of Theorem \ref{th-sl2r} are
somewhat symmetric. As we shall see,  $Q_1$ contains $\{(a,b,c):
a\ge b \ge c \} $ and $Q_2$ contains $\{(a,b,c): a \le b-c
 \} $. Case (3) is of a completely different nature and
it is not even entirely clear, a priori, that it occurs at all. In
the forthcoming work \cite{cgsc08}, we show that $Q_1\cup Q_2$ is
a dense open set in $Q$ and that $S_0$ is an hypersurface
separating $Q_1$ from $Q_2$. This however requires different
techniques than those used in this paper.
\end{remark}

\begin{remark}
In case (3), let $d(t)$ be the metric on $M=\mbox{SL}(2,\mathbb
R)$ induced by $g(t)$. Observe that there are no factors $r(t)$
such that $(M, r(t)d(t))$ converges uniformly to a metric
structure on $M$. A meaningful scaling might be to consider $(M,
\frac{A_0}{A(t)} g(t))$ for which two components converge and the
third goes to zero (potentially, a dimensional collapse but
curvatures blow up).
\end{remark}

\begin{remark}
For the forward normalized Ricci flow (\ref{normRF}), Isenberg and
Jackson \cite{ij92} show that the solution exists for all time and
presents a pancake degeneracy.
\end{remark}

For the proof of Theorem \ref{th-sl2r}, we recall that the
sectional curvatures are
\begin{align*}
K(f_2 \wedge f_3)&=\frac{1}{ABC}(-3A^2+B^2+C^2-2BC-2AC-2AB),\\
K(f_3 \wedge f_1)&=\frac{1}{ABC}(-3B^2+A^2+C^2+2BC+2AC-2AB),\\
K(f_1 \wedge f_2)&=\frac{1}{ABC}(-3C^2+A^2+B^2+2BC-2AC+2AB).
\end{align*}

Therefore, writing
$$g=  Af^1\otimes f^1+Bf^2\otimes f^2+Cf^3\otimes f^3$$
for the solution of the positive normalized  Ricci flow with
initial data $g_0$, $A,B,C$ (with $ABC=4$) satisfy the equations

\begin{equation}\label{pdesl2}
\left \{
\begin{aligned}
\frac{dA}{dt}=&-\frac{2}{3}[-A^2(2A+B+C)+A(B-C)^2],\\
\frac{dB}{dt}=&-\frac{2}{3}[-B^2(2B+A-C)+B(A+C)^2],\\
\frac{dC}{dt}=&-\frac{2}{3}[-C^2(2C+A-B)+C(A+B)^2].
\end{aligned}
\right .
\end{equation}

Without loss of generality we may assume that $B_0 \geq C_0$.
Looking at the evolution equation of $B-C$, it follows that $B
\geq C$ as long as a solution exists.

Since that $ABC=A_0B_0C_0=4$. We have $$\ddt C=-\frac{2}{3}
[BC^2+B^2C-2C^3+ABC-AC^2+ABC+A^2C]\leq -\frac{2}{3},$$  $C$ is
decreasing and the solution can only exist up to some finite time
$T_+<\infty$.

\begin{lemma}Assume that $B_0 \geq C_0$, we have
$$\lim_{T_+} C=0.$$
\end{lemma}
\begin{proof}Observe that $AB$ is increasing because that
$$\ddt \ln (AB)=-\frac{2}{3}[-B^2-A^2-2AB-BC+AC+2C^2]>0.$$

Assume that $\lim_{T_+} C(t)=\eta>0$. Then $AB<\frac{4}{\eta}$ and
$B\ge C>\eta$. So $A<\frac{4}{\eta^2}$ and we must have
$\lim_{T_+} A(t)=0$ and $\lim_{T_+} B(t)=\infty$ (because $AB$ is
increasing and bounded from above, it is easy to see those two
conditions are equivalent). Hence we have
$$\frac{dA}{dt} \sim -\frac{2}{3}AB^2,
\;\; \frac{dB}{dt} \sim \frac{4}{3}B^3,\;\; \frac{dC}{dt} \sim
-\frac{2}{3}CB^2.$$ So we have $B(t)^{-2}\sim  \frac83(T_+-t)$,
but this contradicts $\lim_{T_+} C(t)=\eta>0$.
\end{proof}

\begin{lemma} Assume $B_0\ge C_0$. If there exists a time $t_0$
such that $A(t_0)\ge B(t_0)$ then $$A\sim
\frac{\sqrt{6}}{4}(T_+-t)^{-1/2},\;\; B\sim
\eta_1(T_+-t)^{1/4},\;\; C\sim \eta_2(T_+-t)^{1/4}.$$
\end{lemma}
\begin{proof}
As $ABC=4$, we have $\lim_{T_+} AB=\infty$. Moreover,
\begin{equation}\label{lnab}
\ddt \ln (A/B)=2(A+B)(A+C-B).
\end{equation}

If there exist a time $t_0$ such that $A(t_0)\geq B(t_0)$, then
$A\ge B$ on $[t_0, T_+)$. Similarly, the condition $A>2B$ is
preserved by the flow. Assuming that $A>2B$, we have
$$\ddt B \sim -\frac23 B(A-2B)(A+B)<0,$$ hence $\lim_{T_+}
B(t)=B(T_+)<\infty$. So $\lim_{T_+} A(t)=\infty$ and the system
\ref{pdesl2} yields
\begin{equation}
\left \{
\begin{aligned}
\frac{dA}{dt} \sim &\frac43 A^3,\\
\frac{dB}{dt} \sim &-\frac{2}{3}BA^2,\\
\frac{dC}{dt} \sim &-\frac{2}{3}CA^2.
\end{aligned}
\right .
\end{equation}
This give the desired asymptotics.

We now need to rule out the case when $B\leq A\leq 2B$ for all
$t$. In that case we have $$\frac43 A^3\le \frac{dA}{dt} \le
2A^3,$$ this implies $\int_0^{T_+} A^2 =\infty$. Further, by
(\ref{lnab}), $\int_0^{T_+} (A^2-B^2) <\infty$. Since $A/B$ is
nondecreasing, there exists a constant $\eta$, such that $A-B>\eta
A$. Thus $\int_0^{T_+} A^2 <\infty$. This is a contradiction.
\end{proof}

\begin{lemma} Assume $B_0\ge C_0$. If there exists a time $t_0$
such that $A(t_0) \le B(t_0)-C(t_0)$ then $$A\sim
\eta_1(T_+-t)^{1/4},\;\; B\sim
\frac{\sqrt{6}}{4}(T_+-t)^{-1/2},\;\; C\sim \eta_2(T_+-t)^{1/4}.$$
\end{lemma}

\begin{proof}
We have
\begin{equation}\label{}
\ddt (B-A-C)=-\frac23(2A^3+2C^3-2B^3+2A^2B+2ABC-2AB^2-2A^2C).
\end{equation}
Hence  the condition $B-A\geq C$ is preserved by the flow. It
follows from the flow equation (\ref{pdesl2}) and (\ref{lnab})
that both $B$ and $B/A$ are increasing. So we have $\lim_{T_+}
B=\infty$. If $\lim_{T_+} B/A<\infty$ then, since there exist a
$\eta$ such that $B-A-C>\eta B$, (\ref{lnab}) yields that
$\int_0^{T_+} B^2 <\infty$. The evolution equation (\ref{pdesl2})
shows that this contradicts $\lim_{T_+} B(t)=\infty$. Hence we
must have $\lim_{T_+} B/A=\infty$, and (\ref{pdesl2}) gives
\begin{equation}\label{}
\left \{
\begin{aligned}
\frac{dA}{dt} \sim &-\frac{2}{3}AB^2,\\
\frac{dB}{dt} \sim &\frac{4}{3}B^3,\\
\frac{dC}{dt} \sim &-\frac{2}{3}CB^2.
\end{aligned}
\right .
\end{equation}
This proves the desired result.
\end{proof}

Now the only case left is when $A<B<A+C$ for all $t\in [0,T_+)$.
In this case since $\lim_{T_+} C=0$, we have $\lim_{T_+} (B-A)=0$,
and the flow equation (\ref{pdesl2}) yields that
\begin{equation}\label{}
\left \{
\begin{aligned}
\frac{dA}{dt} \sim &\frac{4}{3}A^3,\\
\frac{dB}{dt} \sim &\frac{4}{3}B^3,\\
\frac{dC}{dt} \sim &-\frac{8}{3}CB^2.
\end{aligned}
\right .
\end{equation} So we arrive at
$$A\sim
\frac{\sqrt{6}}4 (T_+-t)^{-1/2},\;\; B\sim \frac{\sqrt{6}}4
(T_+-t)^{-1/2},\;\; C\sim \frac{32}{3} (T_+-t).$$

This together with the two previous lemmas concludes the proof of
Theorem \ref{th-sl2r}. As noted in the remark following the
theorem, it is not clear from the proof itself that the third case
does indeed occur. In \cite{cgsc08}, we show that there is a
smooth hypersurface of initial condition which is preserved by the
flow, which exactly corresponds to the asymptotic behavior
described in the third case $A<B<A+C$ for all $t\in [0, T_+)$.

\bibliographystyle{halpha}
\bibliography{bio}
\end{document}